\documentclass[a4paper]{article}

 \usepackage{multirow}
\usepackage{amsmath,amssymb,amsfonts}
\usepackage{bm}
\usepackage{graphicx}
\usepackage{accents}
\usepackage{mathtools}

\usepackage{float}
\usepackage{caption}
\usepackage{subcaption}
\usepackage{amsmath}
\usepackage{relsize}
\usepackage{bbm}
\usepackage{epsfig}
\usepackage{mathtools}
\usepackage{algorithm}
\usepackage{algpseudocode}

\input{epsf}

\newtheorem{remark}{Remark}

\def\bkR{{\rm I\kern-.17 em R}}

\def\bkN{{\rm I\kern-0.07em \mbox{\vrule height0.7em
            width0.07em depth0em}\kern-.18em N}}

\begin{document}

\title{Improved PDE Models for Image Restoration  through Backpropagation}

\author{S\'{\i}lvia  Barbeiro\thanks{email: silvia@mat.uc.pt},\quad
Diogo Lobo
\thanks{diogo.lobo@mat.uc.pt}  \\
{\small CMUC, Department of Mathematics, University of Coimbra}\\
 {\small Apartado 3008, EC Santa Cruz, 3001 - 501 Coimbra, Portugal}} 



\date{}

\maketitle

\begin{abstract}

In this paper we focus on learning optimized partial differential equation (PDE) models   for image filtering. 
In this approach, the grey-scaled images are represented by a
vector field of two real-valued functions and the image restoration problem is modelled
by an evolutionary process such that the restored image at any time satisfies an initial-boundary-value problem of cross-diffusion with reaction type. The coupled evolution of the two components of the image is determined by a nondiagonal matrix  that depends on those components.
 A critical question when designing a good-performing filter lies in the selection of the optimal coefficients and influence functions which define the cross-diffusion matrix. 
We propose the use of deep learning techniques in order to optimize the parameters of the model. 
In particular, we use a back propagation technique in order to minimize a cost function related to the quality of the denoising processe, while we ensure stability during the learning procedure. 
Consequently, we obtain improved image restoration models with solid mathematical foundations. The learning framework and resulting models are presented along with related numerical results and image comparisons.
\end{abstract}

\maketitle

\section{Introduction}

Nonlinear diffusion processes are well-known and widely  used for image noise removal. Roughly speaking, the  idea 
is to combine an effective noise reduction by diffusion with the preservation of the edges and other 
important image features. Amongst the better known approaches are those  where the diffusion coefficient depends on the gradient with an inverse proportion, \cite{PerMal90,Wei97}. 

The use of  nonlinear complex diffusion filters (NCDF) where  the image is represented by a complex function and  the process of filtering is governed by a diffusion equation with a complex-valued diffusion coefficient, was investigated in \cite{GilSocZee04}. Those filters, which bring the advantage of using the imaginary part of the solution as an edge detector avoiding the computation of the gradient to control the diffusion coefficient,  can be successfully applied for denoising in particular for medical imaging despeckling \cite{BerEtAl10,SalFer07}.  

A complex diffusion equation can be written as a cross-diffusion system for the real and imaginary parts. 
The use of general cross-diffusion systems for image processing, which encompasses the complex-diffusion equations as particular cases, was investigated in \cite{AraEtAl17b}.
The theoretical studies of the correspondent  initial valued boundary problems were presented  in \cite{AraEtAl17a} and \cite{AraEtAl17b} for the  linear and nonlinear cases, respectively, where  well-posedness and some scale-space representation properties were derived  under hypothesis on the  cross-diffusion coefficient matrix. However, besides these conditions, there is not much insight on the form that these coefficients should take in practical applications. In \cite{AraEtAl17c}, the same authors applied the cross-diffusion model to several examples related to the problem of reducing the speckle noise in optical coherence tomography (OCT) images. The conclusion is that the performance of the model is greatly influenced by the choice of the cross-diffusion coefficient matrix. 

All those methods based on  diffusion filtering are highly dependent on some crucial parameters  ({\it e.g.}  \cite{AraEtAl17b,TsiPet13}). 
The parameters that lead to the most effective methods vary depending on the acquisition methods, the nature of the images and the associated noise profile. Furthermore, diffusion filters require the specification of a stopping time. To circumvent this issue,  a reaction term can be introduced which keeps the steady-state solution close to the original image \cite{Wei97}.

The  challenge  lies in the formulation of robust models with optimal parameters corresponding to each 
architecture.  The use of neural networks for image processing tasks  has set the pace of current research in this area. However, the general intrinsic ``black box" nature of such algorithms is a drawback, particularly in the context of medical applications. 
In \cite{ChePoc17} the authors use a learning framework inspired in nonlinear reaction diffusion processes where the filters and the influence functions are simultaneously learned from training data. The results are auspicious. There are, however, some important questions remaining: the numerical solution is not related with the original reaction diffusion problem, and more importantly, the qualitative properties of the computed solution are not derived. Moreover, the method does not cover promising models as the one based on cross-diffusion. 
 
We propose the use of deep learning techniques such as backpropagation \cite{DL16} and Adam algorithm \cite{Adam14} to optimize the parameters of a system of PDEs for a specific task. In this work, we focus on nonlinear cross-diffusion with reaction filters in order to optimize the reaction parameter and the influence functions of the cross-diffusion matrix for the removal of gaussian noise. 

The optimization process we suggest includes the stability of the inherent numerical method, which is crucial for obtaining successful results.  
We drive our attention first to synthetic images in order to develop our learning strategy. The methodology consists in 
the parametrization of the functions that govern the cross-diffusion system, the computation of the gradient of a cost function with respect to these parameters and the application of the backpropagation technique, widely used in neural network supervised learning, to learn the optimal parameters. Experiments in real noisy images show that the training process improves significantly the performance of the filters. 

Our contributions can be summarized as follows:
\begin{itemize}
	\item we derive the stability conditions for a broad class of cross-diffusion models;
	\item we build a recurrent neural network equivalent to a cross-diffusion system where its neurons are the points of the spatial discretization and the weights are the non-linear diffusion coefficients;
	\item we obtain improved stable cross-diffusion models for image denoising through backpropagation.
\end{itemize}
 
The article is organized as follows. In Section 2 the fully discrete cross-diffusion with reaction model is presented. Next, the stability conditions are derived. Section 4 is devoted to the learning model. We present the numerical experiments in Section 5. We end the paper with a section dedicated to conclusions.  

\section{Cross-diffusion reaction model }

In this section we present a fully discrete cross-diffusion reaction model   for image restoration. 
The image is represented by a two-component vector field,  ${\bf w}=(u,v)^\top$, and the restoration process is governed by the nonlinear cross-diffusion reaction system\begin{equation} \label{cross-diffusion}
\begin{cases}
u_t=\nabla \cdot(d_{1}({\bf w})\nabla u + d_{2}({\bf w})\nabla v) - \lambda (u-u^0)\text{ in }\Omega\times\mathbb{R}^+,\\
v_t=\nabla \cdot(d_{3}({\bf w})\nabla u + d_{4}({\bf w})\nabla v)\text{ in }\Omega\times\mathbb{R}^+,\\
u({\bf x},0)=u^0({\bf x}),\text{ }v({\bf x},0)=v^0({\bf x})\text{ in } \Omega,\\
u_\eta =0,\text{ }v_\eta=0\text{ on } \Gamma\times\mathbb{R}^+,
\end{cases}
\end{equation}
where  $\Omega=(a_{1},b_{1})\times (a_{2},b_{2})\subset\mathbb{R}^{2}$ is the domain of interest, $u^0$ and $v^0$ are the given initial conditions for $u$ and $v$ and $\eta$ denotes the outward normal vector to the boundary  $\Gamma=\partial \Omega$.   
The cross-diffusion matrix  of the model is given by
\begin{equation}\label{cross-diffusion matrix}
D({\bf w})=\begin{pmatrix}d_{1}({\bf w})&d_{2}({\bf w})\\d_{3}({\bf w})&d_{4}({\bf w})\end{pmatrix}.
\end{equation}
$\lambda$ is a time dependent non-negative parameter.

 The domain $\overline{\Omega}=\Omega\cup\Gamma$ is discretized by the points ${\bf x}_{{\bf j}}=(x_{j_{1}},x_{j_{2}})$, where
\[
x_{j_{1}}=a_{1}+h_{1}j_{1},\; x_{j_{2}}=a_{2}+h_{2}j_{2},\quad j_{k}=0,1,\ldots,N_{k},\quad h_{k}=\frac{b_{k}-a_{k}}{N_{k}}, \; k=1,2,
\]
for two given integers $N_{1}, N_{2}\geq 1$, ${\bf j}=(j_{1},j_{2})$ and ${\bf h}=(h_{1},h_{2})$ . This spatial mesh on $\overline{\Omega}$ is denoted by $\overline{\Omega}_{\bf h}$ and  $\Gamma_{\bf h}=\Gamma \cap \overline{\Omega}_{\bf h}$.  
Points halfway between two adjacent grid points are denoted by 
${\bf x}_{{\bf j}\pm (1/2){\bf e}_k}={\bf x}_{\bf j} \pm  \frac{h_k}{2}{\bf e}_{k}$, $k=1,2$, where $\{{\bf e}_{1},{\bf e}_{2}\}$ 
is the canonical basis, that is, ${\bf e}_{k}$ is the standard basis unit vector in the $k$th direction. 

For the discretization in time, we consider a mesh with time step $\Delta t$,
$0=t^0<t^1<t^2<\ldots$, where $t^{m+1}-t^m=\Delta t$. 

We denote by $Z_{\bf j}^m$ the value of a mesh function $Z$ at the point $({\bf x}_{{\bf j}}, t^m)$. 
For the formulation of the finite difference approximations, we use the centered finite difference quotients in the $k$th spatial direction, for $k=1,2$,
\begin{eqnarray*}
 \delta_{k}Z_{\bf j} = \frac{Z_{{\bf j}+(1/2){\bf e}_{k}}-Z_{{\bf j}-(1/2){\bf e}_{k}}}{h_{k}}, \quad \delta_{k}Z_{{\bf j}+(1/2){\bf e}_{k}} = \frac{Z_{{\bf j}+{\bf e}_{k}}-Z_{\bf j}}{h_{k}}.\label{sm2a}
\end{eqnarray*}

In order to formulate the discrete cross-diffusion restoration problem, let $u^0:\Omega_{\bf h}\to\mathbb{R}$ be a discrete real-valued function standing for the grey level values on $\Omega_{\bf h}$ of the noisy image to be restored. From $u^0$, an initial distribution ${\bf W}^{0}=(U^0,V^0)$, for the cross-diffusion, is required. This is given by two real-valued functions $U^0,V^0: \Omega_{\bf h}\rightarrow\mathbb{R}$, that can be selected following different criteria. A simple choice, that will be considered in the experiments we present later in this paper, consists of taking $U^0({\bf x}_{{\bf j}})=u^0({\bf x}_{{\bf j}})$ and $V^0({\bf x}_{{\bf j}})=0$, ${\bf x}_{{\bf j}} \in \overline \Omega_{\bf h}$. A more detailed discussion on the initial data can be seen in \cite{AraEtAl17a,AraEtAl17b}.

Let ${\bf W}^{m}_{\bf j}=(U_{{\bf j}}^{m},V_{{\bf j}}^{m})^\top$, such that ${\bf x}_{\bf j} \in  \overline \Omega_{\bf h}$. Given the initial solution ${\bf W}^{0}_{\bf j}=(U^0_{\bf j},V^0_{\bf j})$,  the numerical solution of (\ref{cross-diffusion}) at the time $t^{m+1}$ is obtained considering the following finite difference scheme  
\begin{eqnarray}
 \frac{U_{\bf j}^{m+1}-U_{\bf j}^{m}}{\Delta t} & = &\sum_{k=1}^{2}\delta_{k}\left(d_{1}({\bf W}^{m})_{\bf j}\delta_{k}U_{{\bf j}}^{m+\theta}+d_{2}({\bf W}^{m})_{\bf j}\delta_{k}V_{{\bf j}}^{m+\theta}\right)\nonumber\\
&&-\lambda^{m+\theta} (U_{\bf j}^{m+\theta}-U_{\bf j}^0) ,\label{cross_schemeEq1}\\
\displaystyle  \frac{V_{\bf j}^{m+1}-V_{\bf j}^{m}}{\Delta t} & = &\sum_{k=1}^{2}\delta_{k}\left(d_{3}({\bf W}^{m})_{\bf j}\delta_{k}U_{{\bf j}}^{m+\theta}+d_{4}({\bf W}^{m})_{\bf j}\delta_{k}V_{{\bf j}}^{m+\theta}\right),\label{cross_schemeEq2}
\end{eqnarray}
where
$$
D({\bf W}^m)_{{\bf j}+ (1/2){\bf e}_{k}}=\frac{D({\bf W}^m_{{\bf j}})+ D({\bf W}^m_{{\bf j}+{\bf e}_{k}})}{2}
$$
and $\theta \in \{0,1\}$. 
If $\theta=0$ then (\ref{cross_schemeEq1}) is an explicit method. If $\theta=1$ then method (\ref{cross_schemeEq1}) is semi-implicit  and its solution is obtained by solving a system of linear  equations. In the next section we will derive the stability properties of both, explicit and semi-implicit, methods. Later,  in sections \ref{Learning procedure} and \ref{Deriving the gradients}, we will discuss their role in our learning strategy: to optimize the cross-diffusion model
we will use the explicit scheme ; to denoise the images we will use the semi-implicit discretization of the optimized model.

Equations  (\ref{cross_schemeEq1}) and (\ref{cross_schemeEq2}), for ${\bf j}$ such that 
 ${\bf x}_{{\bf j}} \in \Gamma_{{\bf h}}$,  are defined using  points of the form ${\bf x}_{{\bf j}}+{\bf e}_{k}$ and 
${\bf x}_{{\bf j}}-{\bf e}_{k}$, which don't belong to $\Omega_{\bf h}$ and are not yet defined. For those points ${\bf x}_{{\bf j}}+{\bf e}_{k}$ and
${\bf x}_{{\bf j}}-{\bf e}_{k}$, we consider the solution in ${\bf x}_{{\bf j}}-{\bf e}_{k}$ and ${\bf x}_{{\bf j}}+{\bf e}_{k}$, respectively. This corresponds to the usual discretization of the homogeneous Neumann boundary conditions on $\Gamma_{{\bf h}}$ with central finite differences.

\section{Stability of the numerical scheme}\label{section_stability}

We will now investigate the stability of the finite difference scheme  (\ref{cross_schemeEq1})-(\ref{cross_schemeEq2}). The approach generalizes the strategies  in \cite{AraBarSer12} used to study the stability of the particular case of  complex-diffusion equations written as a cross-diffusion system.  

For each $x_{\bf j}=(x_{j_1},x_{j_2})\in \bar{\Omega}_{\bf h}$, we define the rectangle $\square_{\bf j}=(x_{j_1},x_{j_1+1})\times(x_{j_2},x_{j_2+1})$ and denote by $|\square_{\bf j}|$ the measure of $\square_{\bf j}$.   We consider the discrete  $L^2$ inner products 
\begin{eqnarray*}\label{inner_product1_2D}
({U},{ V})_h&=&\sum_{\square_{\bf j}\subset \Omega} \frac{|\square_{\bf j}|}{4} \left({ U}_{j_1,j_2}{ V}_{j_1,j_2} +{ U}_{j_1+1,j_2}{{ V}}_{j_1+1,j_2}\right.\\
&&\qquad \left.+{U}_{j_1,j_2+1}{V}_{j_1,j_2+1}+{ U}_{j_1+1,j_2+1}{{ V}}_{j_1+1,j_2+1}\right),
\end{eqnarray*}
\begin{eqnarray*}
({U},{V})_{h_1^*}&=&\sum_{\square_{\bf j}\subset \Omega} \frac{|\square_{\bf j}|}{2} \left({U}_{j_1+1/2,j_2}{{V}}_{j_1+1/2,j_2}
+{ U}_{j_1+1/2,j_2+1}{ V}_{j_1+1/2,j_2+1}\right),
\end{eqnarray*}
and
\begin{eqnarray*}
({ U},{ V})_{h_2^*}&=&\sum_{\square_{\bf j}\subset \Omega} \frac{|\square_{\bf j}|}{2} \left({ U}_{j_1,j_2+1/2}{{ V}}_{j_1,j_2+1/2} +{ U}_{j_1+1,j_2+1/2}{{ V}}_{j_1+1,j_2+1/2}\right).\end{eqnarray*}
Their correspondent norms are denoted by $\|.\|_h$, $\|.\|_{h_1^*}$ and  $\|.\|_{h_2^*}$, respectively. For ${\bf W}=({U},{ V})^\top$ we define $\|{\bf W}\|_h^2=\|{ U}\|_h^2+\|{ V}\|_h^2$.

To simplify the notation and where it is clear from the context, we  write in what follows $d_\ell$ instead of $d_\ell({\bf W}^m),$ or $d_\ell({\bf W}^m)_{\bf j}$, for $\ell=1,2,3,4.$

Multiplying both members of equations (\ref{cross_schemeEq1}) and (\ref{cross_schemeEq2}) by $U^{m+\theta}$ and $V^{m+\theta}$, respectively, according to the discrete inner product $(\cdot,\cdot)_h$, and using summation by parts we get
\begin{equation}
\begin{aligned}
&\left( \frac{{ U}^{m+1}-{ U}^{m}}{\Delta t}, { U}^{m+\theta} \right) _h +\left( \frac{{ V}^{m+1}-{ V}^{m}}{\Delta t}, { V}^{m+\theta} \right) _h+ \sum_{k=1}^2 \Big(
||(d_1)^\frac{1}{2}\delta_k U^{m+\theta}||_{h_k^*}^2 \nonumber\\
&+
||(d_4)^\frac{1}{2}\delta_k V^{m+\theta}||_{h_k^*}^2+(d_2\delta_k V^{m+\theta}, \delta_k U^{m+\theta})_{h_k^*}
+
(d_3\delta_k U^{m+\theta}, \delta_k V^{m+\theta})_{h_k^*}\Big)\nonumber\\
&=-\lambda^{m+\theta}(U^{m+\theta}-U^0,U^{m+\theta})_h.
\end{aligned}
\end{equation}

We can write
\begin{equation*}
{\bf W}^{m+\theta}=\frac{{\bf W}^{m+1}+{\bf W}^{m}}{2}+(\theta-\frac{1}{2})\Delta t\frac{{\bf W}^{m+1}-{\bf W}^{m}}{\Delta t}
\end{equation*}
and then
\begin{equation}
\label{stab1}
\begin{aligned}
&\frac{||{\bf W}^{m+1}||_h^2-||{\bf W}^{m}||_h^2}{2\Delta t}+(\theta-\frac{1}{2})\Delta t \bigg|\bigg|\frac{{\bf W}^{m+1}-{\bf W}^{m}}{\Delta t}\bigg|\bigg|_h^2\\
&+ \sum_{k=1}^2 \Big(
||(d_1)^\frac{1}{2}\delta_k U^{m+\theta}||_{h_k^*}^2 
+
||(d_4)^\frac{1}{2}\delta_k V^{m+\theta}||_{h_k^*}^2+(d_2\delta_k V^{m+\theta}, \delta_k U^{m+\theta})_{h_k^*}\\
&+
(d_3\delta_k U^{m+\theta}, \delta_k V^{m+\theta})_{h_k^*}\Big)\leq \lambda^{m+\theta}\|U^{m+\theta}-U^0\|_h \|U^{m+\theta}\|_h.
\end{aligned}
\end{equation}
For any $\epsilon>0$ we have that 
\begin{equation}\label{ineq_epsilon}
\begin{aligned}
\lambda^{m+\theta}\|U^{m+\theta}-U^0\|_h \|U^{m+\theta}\|_h &\leq \lambda^{m+\theta}(  \|U^{m+\theta}\|_h+\|U^0\|_h)   \|U^{m+\theta}\|_h\\
 &\leq  \lambda^{m+\theta}\|U^{m+\theta}\|_h^2+ \frac{(\lambda^{m+\theta})^2}{4 \epsilon} \|U^0\|_h^2 + \epsilon  \|U^{m+\theta}\|_h^2.
\end{aligned}
\end{equation}

\subsection{Stability of the semi-implicit scheme}
We start by considering the case where $\theta=1$. With the assumption that 
\begin{equation}\label{stability_inequality_semi}
\begin{aligned}
& \sum_{k=1}^2 \Bigl( ||(d_1)^\frac{1}{2}\delta_k U^{m+1}||_{h_k^*}^2 +||(d_4)^\frac{1}{2}\delta_k V^{m+1}||_{h_k^*}^2\\
&+(d_2\delta_k V^{m+1}, \delta_k U^{m+1})_{h_k^*}+(d_3\delta_k U^{m+1}, \delta_k V^{m+1})_{h_k^*}\Bigr) >0 ,
\end{aligned}
\end{equation}
which corresponds to the natural assumption of the diffusion matrix $D$ to be uniformly positive definite, from 
(\ref{stab1}) and (\ref{ineq_epsilon}) we get
\begin{equation*}
(1- 2\Delta t \lambda^{m+1}- 2\Delta t \epsilon)  ||{\bf W}^{m+1}||_h^2 \leq   ||{\bf W}^{m}||_h^2  + 2\Delta t  \frac{(\lambda^{m+1})^2}{4 \epsilon}  \|U^0\|_h^2.
\end{equation*}
Assuming that 
\begin{equation}\label{stability_lambda}
\lambda^{m+1} \leq \lambda_{\text{max}}, \quad
 0<\zeta<1- 2\Delta t \lambda^{m+1}- 2\Delta t \epsilon, \quad m=1,2,\ldots
\end{equation}
for some $ \lambda_{\text{max}}, \zeta \in \bkR^+$ where $\epsilon$ is a constant arbitrarily chosen, we get
\begin{equation*}
  ||{\bf W}^{m+1}||_h^2 \leq  (1+ 2\Delta t (\lambda+\epsilon) \zeta^{-1} )||{\bf W}^{m}||_h^2  + \Delta t  \frac{ \lambda_{\text{max}}^2}{2 \epsilon \zeta}  \|U^0\|_h^2.
\end{equation*}
If (\ref{stability_lambda}) holds, using the Duhamel's principle (\cite{ChanShen87}, Lemma 4.1 in Appendix B) we get
 \begin{equation}
\begin{aligned}
||{\bf W}^{m+1}||_h^2\leq e^{  (1+ 2( \lambda_{\text{max}}+\epsilon) \zeta^{-1} ) t^{m+1}} \Bigl(1+t^{m+1}  \frac{ \lambda_{\text{max}}^2}{2 \epsilon \zeta}  \Bigr) ||{\bf W}^{0}||_h^2\nonumber
\end{aligned}
\end{equation}
and we conclude that the method is stable.

\subsection{Stability of the explicit scheme}
We now consider the case where $\theta=0$.
Applying the $||\cdot||_h$ norm on both sides of  equations  (\ref{cross_schemeEq1}) and (\ref{cross_schemeEq2}), and using the inequalities $(a\pm b)^2\leq 2a^2+2b^2$ and $(a\pm b)^2\leq (1+\eta) a^2+ (1+\eta^{-1})b^2$ for $\eta>0$, we obtain, respectively,
\begin{equation}\label{ineqstab1}
\begin{aligned}
\bigg|\bigg|\frac{U^{m+1}-U^{m}}{\Delta t}\bigg|\bigg|_h^2\leq &  \sum_{k=1}^2\frac{8(1+\eta)}{h_k^2}\Big(||d_1\delta_k U^m||_{h_k^*}^2+||d_2\delta_k V^m||_{h_k^*}^2\\
&+2(d_1\delta_k U^m,d_2 \delta_k V^m)_{h_k^*}\Big)\\
&+2(1+\eta^{-1})(\lambda^m)^2 (\|U^m\|_h^2+\|U^0\|_h^2),
\end{aligned}
\end{equation}
for any $\eta>0$ and
\begin{equation}\label{ineqstab2}
\bigg|\bigg|\frac{V^{m+1}-V^{m}}{\Delta t}\bigg|\bigg|_h^2\leq   \sum_{k=1}^2 \frac{8}{h_k^2}\Big(||d_3\delta_k U^m||_{h_k^*}^2+||d_4\delta_k V^m||_{h_k^*}^2 + 2(d_3\delta_k U^m,d_4 \delta_k V^m)_{h_k^*}\Big).
\end{equation}

Combining the inequalities (\ref{ineq_epsilon}), (\ref{ineqstab1}) and (\ref{ineqstab2}) with (\ref{stab1}) leads to
\begin{equation*}
\begin{aligned}
&\frac{||{\bf W}^{m+1}||_h^2-||{\bf W}^{m}||_h^2}{2\Delta t}+  \sum_{k=1}^2 \Bigl( ||(d_1)^\frac{1}{2}\delta_k U^m||_{h_k^*}^2 +||(d_4)^\frac{1}{2}\delta_k V^m||_{h_k^*}^2\\
&+(d_2\delta_k V^m, \delta_k U^m)_{h_k^*}+(d_3\delta_k U^m, \delta_k V^m)_{h_k^*}\\
&-\frac{4 \Delta t}{ h_k^2}\big((1+\eta_1)(||d_1\delta_k U^m||_{h_k^*}^2+||d_2\delta_k V^m||_{h_k^*}^2+2(d_1\delta_k U^m,d_2 \delta_k V^m)_{h_k^*})\\
&+||d_3\delta_k U^m||_{h_k^*}^2+||d_4\delta_k V^m||_{h_k^*}^2+2(d_3\delta_k U^m,d_4 \delta_k V^m)_{h_k^*}\big) \Bigr)\\ & \leq \Big( \lambda^m+\epsilon + \Delta t(1+\eta^{-1})(\lambda^m)^2\Big)  \|U^m\|_h^2\\
&+ \Big(\frac{(\lambda^m)^2}{4 \epsilon} + \Delta t  (1+\eta^{-1})(\lambda^m)^2\Big)  \|U^0\|_h^2.
\end{aligned}
\end{equation*}

Let $\eta=\epsilon$.
In order to obtain a stable scheme we impose that
\begin{equation}\label{stability_inequality}
\begin{aligned}
& \sum_{k=1}^2 \Bigl( ||(d_1)^\frac{1}{2}\delta_k U^m||_{h_k^*}^2 +||(d_4)^\frac{1}{2}\delta_k V^m||_{h_k^*}^2\\
&+(d_2\delta_k V^m, \delta_k U^m)_{h_k^*}+(d_3\delta_k U^m, \delta_k V^m)_{h_k^*}\\
&-\frac{4 \Delta t}{ h_k^2}\big((1+\epsilon)(||d_1\delta_k U^m||_{h_k^*}^2+||d_2\delta_k V^m||_{h_k^*}^2+2(d_1\delta_k U^m,d_2 \delta_k V^m)_{h_k^*})\\
&+||d_3\delta_k U^m||_{h_k^*}^2+||d_4\delta_k V^m||_{h_k^*}^2+2(d_3\delta_k U^m,d_4 \delta_k V^m)_{h_k^*}\big) \Bigr) \geq 0,
\end{aligned}
\end{equation}
for some $\epsilon>0$. 

If (\ref{stability_inequality}) holds, then
\begin{equation}\label{stability_inequality_2}
\begin{aligned}
||{\bf W}^{m+1}||_h^2 \leq &  \|{\bf W}^m\|_h^2+ 2 \Delta t \Big( \lambda^m+\epsilon + \Delta t(1+\epsilon^{-1})(\lambda^m)^2\Big)  \|U^m\|_h^2 \\
&+  2 \Delta t  \Big(\frac{(\lambda^m)^2}{4 \epsilon} +\Delta t (1+\epsilon^{-1})(\lambda^m)^2\Big)  \|U^0\|_h^2. 
\end{aligned}
\end{equation}
We now take $a_\epsilon=2  \Big( \lambda^m+\epsilon + \Delta t(1+\epsilon^{-1})(\lambda^m)^2\Big) $ and  $b_\epsilon=
 2  \Big(\frac{(\lambda^m)^2}{4 \epsilon} +\Delta t (1+\epsilon^{-1})(\lambda^m)^2\Big)$. From (\ref{stability_inequality_2}) and using the Duhamel's principle  we get
 \begin{equation}
\begin{aligned}
||{\bf W}^{m+1}||_h^2\leq e^{a_\epsilon t^{m+1}} \Bigl(1+t^{m+1} b_\epsilon \Bigr) ||{\bf W}^{0}||_h^2\nonumber
\end{aligned}
\end{equation}
and we conclude that the method is stable.

Note that we can rewrite (\ref{stability_inequality})  as
\begin{equation}\label{stability_inequality2}
\displaystyle \sum_{k=1}^2(\delta_k W)^{m\top}M \delta_k W^m\geq 0,
\end{equation}
where $M$ is a square matrix of dimension $2N_1N_2\times 2N_1N_2$.
In order for (\ref{stability_inequality2}) to hold we require $M$ to be semi-positive definite, that is, the eigenvalues of $M$ are all non-negative. Using the Gershgorin Theorem,  a sufficient condition for the method to be stable  is to impose that the influence functions satisfy, for some $\epsilon>0$,
\begin{equation}
\label{stablity_gresgorovitch}
\begin{aligned}
&d_1-\frac{4\Delta t}{h_k^2}((1+\epsilon)d_1^2+d_3^2)\geq |\frac{1}{2}(d_2+d_3)-\frac{4\Delta t}{h_k^2}( (1+\epsilon) d_1d_2+d_3d_4)|,\\
&d_4-\frac{4\Delta t}{h_k^2}((1+\epsilon)d_2^2+d_4^2)\geq |\frac{1}{2}(d_2+d_3)-\frac{4\Delta t}{h_k^2}((1+\epsilon) d_1d_2+d_3d_4)|.
\end{aligned}
\end{equation}

\section{Learning model}

Our goal is to optimize the parameters of our model using deep learning techniques. 
In order to adapt the cross-diffusion scheme (\ref{cross_schemeEq1})-(\ref{cross_schemeEq2}) to a neural network architecture, we concatenate $U$ and $V$ into a single column vector $w$ with $2N_1N_2$ entries:
\begin{equation*}w=(u,v)^\top,
	\quad u=
	\begin{bmatrix}
		U_{1,1}\\U_{1,2}\\\vdots \\ U_{N_1,N_2}
	\end{bmatrix}, \quad
	v=
	\begin{bmatrix}
		V_{1,1}\\V_{1,2}\\\vdots \\ V_{N_1,N_2}
	\end{bmatrix}.
\end{equation*}

The $m$-th iteration of the scheme (\ref{cross_schemeEq1})-(\ref{cross_schemeEq2})  can now be written in the vector formulation
\begin{eqnarray}\label{learning_scheme}
	\lefteqn{\frac{w^{m+1}-w^m}{\Delta t}}\nonumber\\
	&=&(K_l^xD_L^{x,m}K_{rL}^x+K_l^yD_L^{y,m}K_{rL}^{y}+K_l^xD_R^{x,m}K_{rR}^x+K_l^yD_R^{y,m}K_{rR}^y )w^{m+\theta}\nonumber\\
	&&-\lambda^{m+\theta} (u^{m+\theta}-u^0),
\end{eqnarray}
with $\lambda^{m+\theta}>0$, and where the first difference matrices ($K_{r*}^*$) are block matrices:
\begin{equation*}
	K_{rL}^x=
	\begin{bmatrix}
		k_r^x & 0 \\ 0 & k_r^x
	\end{bmatrix},
	\quad
	K_{rR}^x=
	\begin{bmatrix}
		0 & k_r^x \\ k_r^x & 0 
	\end{bmatrix},
	\quad
	K_{rL}^y=
	\begin{bmatrix}
		k_r^y & 0 \\ 0 & k_r^y
	\end{bmatrix},
	\quad
	K_{rR}^y=
	\begin{bmatrix}
		0 & k_r^y \\ k_r^y & 0
	\end{bmatrix},
\end{equation*}
 $k_r^x$ and $k_r^y$ denote the backward difference operators with respect to $x$ and $y$, respectively.
The  matrices $D_L^{x,m},D_L^{y,m},D_R^{x,m}$ and $D_R^{y,m}$  are  $2N_1N_2\times 2N_1N_2$ diagonal matrices. 
Each entry $(j,j)$ of $D_L^{x,m}$ and $D_L^{y,m}$ depend on the function $d_1$ for $j=1,\ldots,N_1N_2$  and on the function $d_4$ for $j=N_1N_2+1,\ldots,2 N_1N_2$, at the time $t^m$. The entries of $D^{x,m}_R$ and $D^{y,m}_R$ follow the same pattern, replacing $d_1$ and $d_4$ by $d_2$ e $d_3$, respectively.
It remains to define the difference matrices $K_l^x$ and $K_l^y$:
\begin{equation*}
	K_l^x=
	\begin{bmatrix}
		k_l^x & 0 \\ 0 & k_l^x
	\end{bmatrix},
	\quad
	\quad
	K_l^y=
	\begin{bmatrix}
		k_l^y & 0 \\ 0 & k_l^y
	\end{bmatrix},
\end{equation*}
where $k_l^x$ and $k_l^y$ denote the forward difference operators with respect to $x$ and $y$, respectively.

We reiterate that the scheme (\ref{learning_scheme}) is equivalent to the cross-diffusion scheme (\ref{cross_schemeEq1})-(\ref{cross_schemeEq2}).

\subsection{Parameterization of the influence functions}

An important aspect of the methodology of our learning strategy consists in the parametrization of arbitrary influence functions. 
We parameterize the  influence functions of the cross-diffusion matrix (\ref{cross-diffusion matrix}) through gaussian radial basis functions (RBFs) \cite{Buh03}. Here we consider that they only depend on the second component of the vector $w=(u,v)^\top$, {\it i.e.},  they depend on the component which plays the role of edge detector.  Each of these functions has the expression
\begin{equation}
	\label{rbf}
	d_\ell(v)=\sum_{i=1}^{P}\delta_{\ell,i}\phi \bigg(\frac{|v-\mu_i|}{2\nu}\bigg),\quad \phi(z)=e^{-||z||^2},
\end{equation}
where $\mu_i$, for $i=1,...,P$, are equidistant points in the set $A$ of edge-detector range of values, $\nu$ is a positive scalar (the so called scale of the RBF), and $\delta_{\ell,i}$ are, for each $\ell=1,\dots,4$, the $P$ interpolants od $d_\ell$. A good approximation requires a careful balance between $\nu$ and $P$.

In the learning model, we consider in (\ref{learning_scheme}) $\lambda^{m+\theta}=\lambda$, where $\lambda$ is a positive scalar.

\subsection{Constrained optimization} \label{Learning procedure}
In order to optimize the  influence functions of the cross diffusion matrix $d_{1},d_{2},d_{3}$ and $d_{4}$ and also the scalar $\lambda$ which characterizes the reaction term,  a training cost function and a training dataset are required. We will use a set of $B$ gray-scale images which will serve as basis for our training set. 
We fix the time-step $\Delta t$ and a number of steps $M$, which will result in a cross-diffusion stopping time of $T=M\Delta t$. In the learning viewpoint, we will have a $M$-layered convolutional neural network. We aim to minimize the loss function
$L$, 
\begin{equation}
\label{loss_function}
	L(\Theta,w_1^m,\dots,w_B^m,w_1^{nl},\dots,w_B^{nl})=\sum_{i=1}^{B}l(\Theta,w_i^m,w_i^{nl})= \sum_{i=1}^{B}\frac{1}{2}||u_i^m-u_i^{nl}||^2,
\end{equation}
with $w_i^m=(u_i^m,v_i^m)^\top$, $w_i^{nl}=(u_i^{nl},v_i^{nl})^\top$, where
 $w_i^m$ and $w_i^{nl}$ are, respectively, the $m$- th iteration of (\ref{learning_scheme}) of the $i$-th corrupted image and the non-corrupted image, respectively, $B$ is the learning batch size, and $\Theta$ is the set of parameters we want to optimize. As such, the set of parameters $\Theta$ is
\begin{equation*}
	\label{parameters}
	\Theta=\{\lambda,\delta_{1,1},\dots,\delta_{1,P},\delta_{4,1},\dots,\delta_{4,P},\delta_{2,1},\dots,\delta_{2,P},\delta_{3,1},\dots,\delta_{3,P}\},
\end{equation*}
and we must now obtain $\frac{\partial L}{\partial \Theta}(\Theta,w_1^M,\dots,w_B^M,w_1^{nl},\dots,w_B^{nl})$ to solve the minimization problem. 

We need to guarantee that the learning procedure ends up with a stable scheme. For that, we will make use of the stability conditions derived in Section \ref{section_stability}. Although the use of the explicit method (\ref{learning_scheme}) with $\theta=0$ can bring some advantages in terms of computational effort, the stability conditions are much more restrictive when compared to the semi-implicit method (\ref{learning_scheme}) with $\theta=1$. Moreover, since we need to impose the stability conditions over the values of $\Theta$, the nonlinearities in the stability condition (\ref{stablity_gresgorovitch}) corresponding to the explicit scheme carry out some issues that we want to avoid here. For that reason our goal is to learn a semi-implicit scheme for image denoising. 

To obtain a  stable semi-implicit scheme, instead of minimizing the loss function (\ref{loss_function}) we seek the solution of the constrained optimization problem
\begin{subequations}\label{constraint}
	\begin{alignat}{2}
		&\!\min_{\Theta}        &\qquad& L(\Theta,w_1^M,\dots,w_B^M,w_1^{nl},\dots,w_B^{nl})\label{opt_problem}\\
		&\text{s.t.} &      & d_1(x)\geq \tfrac{1}{2}|d_2(x)+d_3(x)|,\quad \forall x\in A,\label{constraint1}\\
		&                  &      & d_4(x)\geq \tfrac{1}{2}|d_2(x)+d_3(x)|,\quad \forall x\in A.\label{constraint2}
	\end{alignat}
\end{subequations}

This constraints ensure that (\ref{stability_inequality_semi}) hold. We first note that (\ref{constraint1})-(\ref{constraint2}) are not in the standard formulation over the values of $\Theta$. As the radial basis functions (\ref{rbf}) are strictly positive, we replace (\ref{constraint1})-(\ref{constraint2}) with
\begin{subequations}
	\begin{alignat}{2}
		&c_{1,i}(\Theta)=\delta_{1,i}-\tfrac{1}{2}(\delta_{2,i}+\delta_{3,i})\geq0,\quad &i=1,\dots,P,\label{dconstraint1}\\
		&c_{2,i}(\Theta)=\delta_{1,i}+\tfrac{1}{2}(\delta_{2,i}+\delta_{3,i})\geq0,\quad &i=1,\dots,P,\label{dconstraint2}\\
		&c_{3,i}(\Theta)=\delta_{4,i}-\tfrac{1}{2}(\delta_{2,i}+\delta_{3,i})\geq0,\quad &i=1,\dots,P,\label{dconstraint3}\\
		&c_{4,i}(\Theta)=\delta_{4,i}+\tfrac{1}{2}(\delta_{2,i}+\delta_{3,i})\geq0,\quad &i=1,\dots,P\label{dconstraint4},
		&\end{alignat}
\end{subequations}
using the fact that (\ref{dconstraint1})-(\ref{dconstraint4}) implies (\ref{constraint1})-(\ref{constraint2}). To solve the inequality constrained problem we define the augmented Lagrangian $\mathcal{L}$ as
\begin{equation}
	\mathcal{L}(\Theta,\mu,\rho)=L(\Theta)+\frac{\rho}{2}\sum_{i=1}^{P}\sum_{\ell=1}^{4}\max\bigg(0,\frac{\mu_{\ell,i}}{\rho}-c_{\ell,i}(\Theta)\bigg)^2,
\end{equation}
which is the usual Powell-Hestenes-Rockafellar function applied to problem (\ref{opt_problem}) with constraints  (\ref{dconstraint1})-(\ref{dconstraint4}). The vector $\mu$ is the set of Lagrange multipliers associated with constraints (\ref{dconstraint1})-(\ref{dconstraint4}) and $\rho$ is the penalty parameter.

We follow the main model algorithm in \cite{Birgin09} with some minor modifications. The learning procedure will iterate through $\mathcal{L}(\Theta_k,\mu_k,\rho_k)$. At the iteration $k$ of the minimization algorithm, the Lagrange multipliers are updated through the formula
\begin{equation}
\label{lagrange_updates}
\mu_{k+1}=\min\{\max\{0, \mu_k-\rho_k c(\Theta_k) \},\bar{\mu}\}
\end{equation}
where $\bar{\mu}>0$ is a multiplier cap, and the penalty parameter is increased or decreased if the infeasibility increases or decreases, respectively. The measure of infeasibility at iteration $k$ is given by the quantity
\begin{equation}
\label{infeasibility}
	I_k=\min\{c(\Theta_k),\mu_k/\rho_k \}. 
\end{equation}
The penalty parameter is updated in the following way: for some predefined $\tau \in (0,1)$ and $\gamma>1$,  if $I_k\leq \tau I_{k-1}$ then $\rho_{k+1}=\rho_k / \gamma$, otherwise $\rho_{k+1}=\gamma  \rho_k$.

\begin{remark}
The stability condition  (\ref{stability_lambda}) is not very limitative and for that reason we don't include it in the formulation of the optimization problem  (\ref{constraint}). Even so, we verify if the condition is satisfied by the optimal solution. In our numerical experiments, presented later in the paper, that was always the case. 
\end{remark}

\subsection{Deriving the gradients}\label{Deriving the gradients}
A central step to solve the optimization problem is to compute the gradients of the loss function with respect to the training parameters.
We note that the direction of maximum growth of $L$ can be obtained via back-propagation:
\begin{equation}
	\label{gradient}
	\begin{aligned}
		\frac{\partial l(\Theta,w^M,w^{nl})}{\partial\Theta}=\sum_{m=1}^{M}\frac{\partial l(\Theta,w^M,w^{nl})}{\partial w^M}\frac{\partial w^M}{\partial w^{M-1}}\cdots\frac{\partial w^m}{\partial \Theta}.
	\end{aligned}
\end{equation}

Although the objective is optimize a semi-implicit cross-diffusion with reaction scheme, we use the explicit scheme  
(\ref{learning_scheme}) with $\theta=0$ for the learning steps. 

We can readily determine
\begin{equation*}
	\frac{\partial l(\Theta,w^M,w^{nl})}{\partial w^M}= \begin{bmatrix}
	(u^{M} -u^{nl})^\top \quad \quad 0\end{bmatrix}.
\end{equation*}

We will obtain $\frac{\partial w^m}{\partial \Theta}$ and $\frac{\partial w^{m+1}}{\partial w^{m}}$ using (\ref{learning_scheme}). We start by noticing that
\begin{equation*}
\frac{\partial w^{m}}{\partial \lambda}=-\Delta t \begin{bmatrix}(u^{m-1}-u^0) \\ 0\end{bmatrix}.
\end{equation*}
We observe that
\begin{equation*}
	K_l^xD_L^{x,m}K_{rL}^xw^m=K_l^x\text{diag}(K_{rL}^xw^{m})g_L^{x,m}
\end{equation*}
where $g_L^{x,m}$ is the vector with the entries of $D_L^{x,m}$. From (\ref{rbf}), 
 we can write $g^{x,m}_L=G^{x,m} \Lambda_{1,4}$, where $\Lambda_ {1,4}$ is a $2P$ sized vector with the parameters  $\delta_{1, i}$, $i=1,\ldots,P$ and  $\delta_{4, i}$,  $i=1,\ldots,P$ of $d_1$ and $d_4$, respectively, and 
\begin{equation*} 
G^{x,m}=\frac{1}{2} \begin{bmatrix}
S^x \Phi ^{x,m} & 0 \\ 0 & S^x\Phi ^{x,m}
\end{bmatrix},
\end{equation*}
where $S^x=|k^x_r|$ and
 with $\Phi ^{x,m}$ being a $N_1N_2\times P$ matrix 
\begin{equation*}
	\Phi ^{x,m}=
	\begin{bmatrix}
		\phi_1(v^m_1)& \cdots & \phi_P(v^m_1)\\
		\vdots &  & \vdots  \\
		\phi_1(v^m_{N_1N_2})& \cdots & \phi_P(v^m_{N_1N_2})\\
			\end{bmatrix}.
\end{equation*}

Analogously, we write $g^{y,m}_L=G^{y,m} \Lambda_{1,4}$, $g^{x,m}_R=G^{x,m} \Lambda_{2,3}$, $g^{y,m}_R=G^{y,m} \Lambda_{2,3}$, where  $S^y=|k^y_r|$, and obtain
\begin{equation*}
\frac{\partial w^m}{\partial \Lambda_{1,4}}= \Delta tK_l^x\text{diag}(K_{rL}^xw^{m-1})G^{x,m-1}+\Delta tK_l^y\text{diag}(K_{rL}^yw^{m-1})G^{y,m-1},
\end{equation*}
\begin{equation*}
\frac{\partial w^m}{\partial \Lambda_{2,3}}= \Delta tK_l^x\text{diag}(K_{rR}^xw^{m-1})G^{x,m-1}+\Delta tK_l^y\text{diag}(K_{rR}^yw^{m-1})G^{y,m-1}.
\end{equation*}
Consequently,
\begin{equation}
	\label{gradient_param}
	\frac{\partial w^m}{\partial \Theta} = 
	\begin{bmatrix} \displaystyle 
	\frac{\partial w^m}{\partial \lambda} \quad \quad \frac{\partial w^m}{\partial \Lambda_{1,4}} \quad \quad \frac{\partial w^m}{\partial \Lambda_{2,3}}
	\end{bmatrix}.
\end{equation}

To derive the expression of $\displaystyle \frac{\partial w^{m+1}}{\partial w^m}$, we consider (\ref{learning_scheme}).  We have 
\begin{equation}
	\label{dws+1_dws}
	\frac{\partial w^{m+1}}{\partial w^m}=I+\frac{\partial w^{m+1}}{\partial w^m}\bigg]_L^x+\frac{\partial w^{m+1}}{\partial  w^m}\bigg]_R^x+\frac{\partial w^{m+1}}{\partial  w^m}\bigg]_L^y+\frac{\partial w^{m+1}}{\partial  w^m}\bigg]_R^y,
\end{equation}
with
\begin{equation*}
	\frac{\partial w^{m+1}}{\partial  w^m}\bigg]_L^x=\Delta tK_l^xD_L^{x,m}K_{rL}^x+\frac{\Delta t}{2}K_l^x\text{diag}(K_{rL}^xw^m)\begin{bmatrix}
		0 & S^xd_1'(v^m) \\ 0 & S^xd_4'(v^m)
	\end{bmatrix}\bigg),
\end{equation*}
\begin{equation*}
	\frac{\partial w^{m+1}}{\partial  w^m}\bigg]_R^x=\Delta tK_l^xD_R^{x,m}K_{rR}^x+\frac{\Delta t}{2}K_l^x\text{diag}(K_{rR}^xw^m)\begin{bmatrix}
		0 & S^xd_2'(v^m) \\ 0 & S^xd_3'(v^m)
	\end{bmatrix}\bigg),
\end{equation*}
\begin{equation*}
	\frac{\partial w^{m+1}}{\partial  w^m}\bigg]_L^y=\Delta tK_l^yD_L^{y,m}K_{rL}^y+\frac{\Delta t}{2}K_l^y\text{diag}(K_{rL}^yw^m)\begin{bmatrix}
		0 & S^yd_1'(v^m) \\ 0 & S^yd_4'(v^m)
	\end{bmatrix}\bigg),
\end{equation*}
\begin{equation*}
	\frac{\partial w^{m+1}}{\partial  w^m}\bigg]_R^y=\Delta tK_l^yD_R^{y,m}K_{rR}^y+\frac{\Delta t}{2}K_l^y\text{diag}(K_{rR}^yw^m)\begin{bmatrix}
		0 & S^xd_2'(v^m) \\ 0 & S^xd_3'(v^m)
	\end{bmatrix}\bigg),
\end{equation*}
where $d_\ell'(v^m)$ is a vector of dimension $N_1N_2$  with entries $[d_\ell'(v^m)]_j=d_\ell'(v_j^m)$, being $v_j^m$ is the second component of the vector $w_j^m=(u_j^m,v_j^m)^\top$. We obtain the derivatives of the influence functions numerically with a centered difference scheme.

\subsection{Back-propagation algorithm}

{\small
\begin{algorithm}
	\caption{Algorithm for NCDF learning}
\begin{algorithmic}[1]
	\Require $B$ gray-scale images, $\sigma>0$ (s.d. of gaussian noise);
	\Require $\triangle t>0$, $(N_1,N_2,M)\in\mathbb{N}^3$, $\Theta_1\in \mathbb{R}^{4P+1}$;
	\Require $K_{\text{max}}>0$, $\bar{\mu}>0$, $\mu_1\in [0,\bar{\mu}]^{4P}$, $\rho_1>0$, $\tau\in ]0,1]$, $\gamma>1$;
	\Require Initialization of Adam parameters. 
	\State Set $k=1$;
	\While{$k\leq K_{\text{max}}$}
		\State Extract $B$ random $N_1\times N_2$ sections from the training set (batch);
		\State Add random noise to the batch;
		\State Compute $\frac{\partial \mathcal{L}}{\partial \Theta_k}$ through backpropagation;
		\State Update $\Theta_k$ via Adam algorithm as in \cite{Adam14};
		\State $\mu_{k+1}=\min\{\max\{0, \mu_k-\rho_k c(\Theta_k) \},\bar{\mu}\}$;
		\If{$I_k\leq \tau I_{k-1}$}
			\State $\rho_{k+1}=\rho_k / \gamma$;
		\Else
			\State $\rho_{k+1}=\gamma  \rho_k$;
		\EndIf
		\State $k=k+1$;
	\EndWhile
	\Ensure $\Theta_{K_{\text{max}}}$
	\label{alg:learning_alg}
\end{algorithmic}
\end{algorithm}}

We now describe how the unknown parameters are learned via back-propagation. We store $B$ gray-scale images to form the basis for our training set. In each iteration, we crop a random $N_1\times N_2$ section of each image and add the desired synthetic random noise to each section (for example noise of Gaussian type with zero mean and standard deviation $\sigma$). This set of freshly cut sections with noise will make up the batch for the iteration.

There is also the need to fix the step time $\Delta t$ and the number of layers $M$ before solving the minimization problem. This is a crucial choice, as $T=M\Delta t$ becomes the diffusion time for which the procedure will be optimized. The functions $d_1,d_2,d_3$ and $d_4$ are initialized as in the nonlinear complex diffusion case \cite{AraEtAl17b}, while $\lambda$ is initialized with the value that provides the best average peak to signal noise ratio (PSNR) for the nonlinear complex diffusion process over the training set.

Given these particular considerations, we present the learning algorithm model (Algorithm 1), which is the Adam algorithm \cite{Adam14} applied to our particular problem.

\section{Experimental results}

\begin{figure}[!htb]
\center
	\minipage{0.9\textwidth}
	\includegraphics[width=\linewidth]{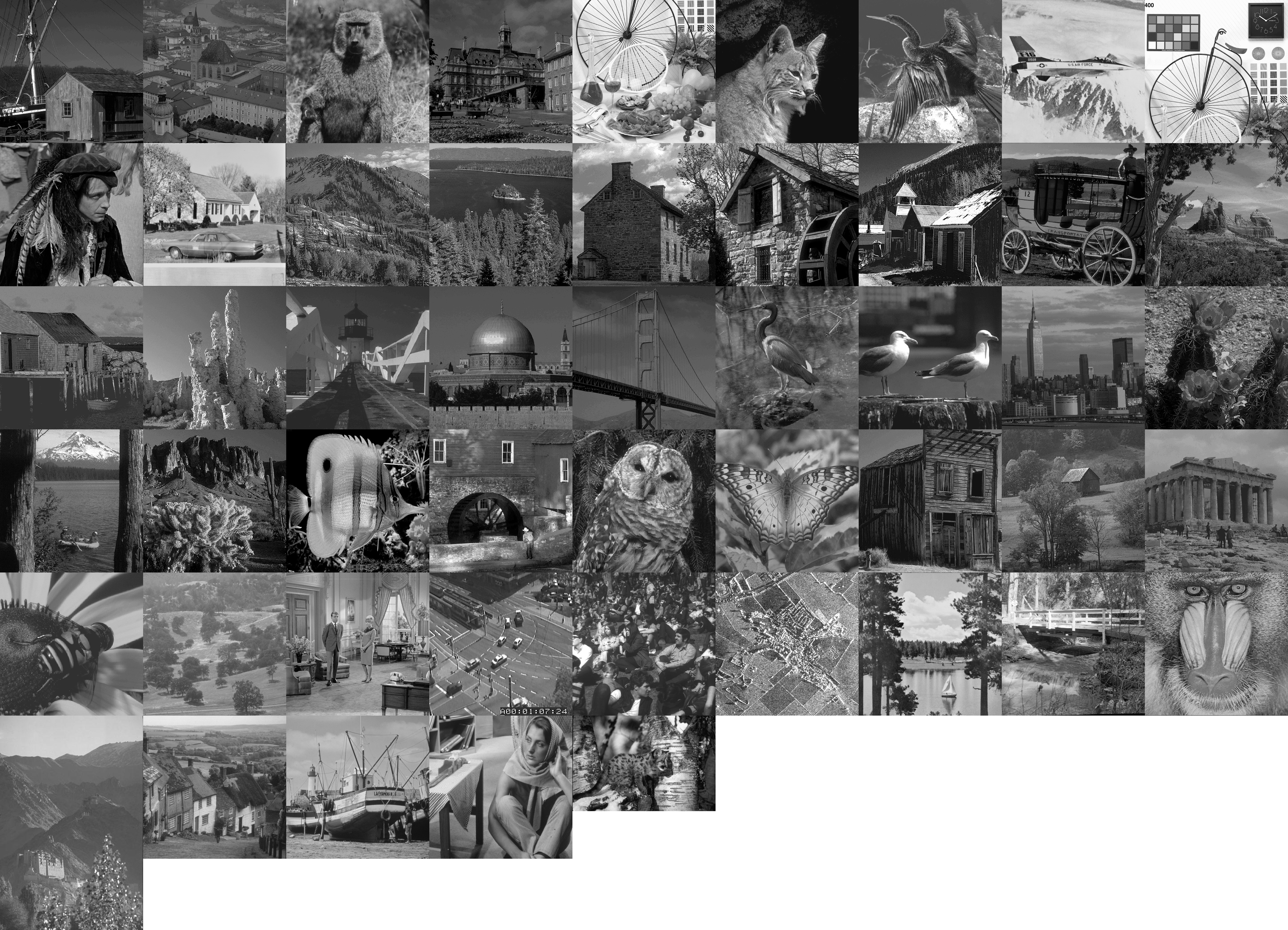}
	\endminipage
	\caption{Images used as training set for the numerical experiments.}
	\label{fig:training_set}
\end{figure}

\begin{figure}[!htb]
\center
	\minipage{0.2\textwidth}
	\includegraphics[width=\linewidth]{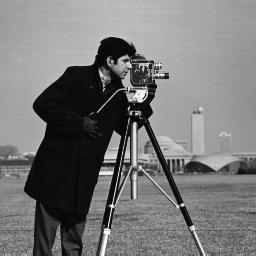}
	\endminipage
	\minipage{0.2\textwidth}
	\includegraphics[width=\linewidth]{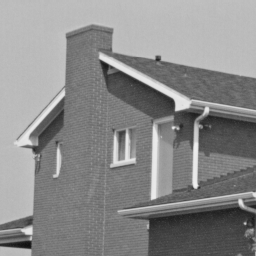}
	\endminipage
	\minipage{0.2\textwidth}
	\includegraphics[width=\linewidth]{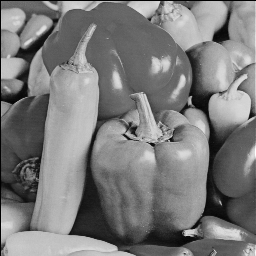}
	\endminipage \\
	
	\minipage{0.2\textwidth}
	\includegraphics[width=\linewidth]{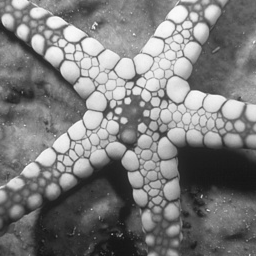}
	\endminipage
	\minipage{0.2\textwidth}
	\includegraphics[width=\linewidth]{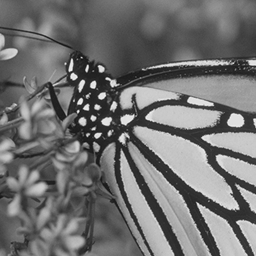}
	\endminipage
	\minipage{0.2\textwidth}
	\includegraphics[width=\linewidth]{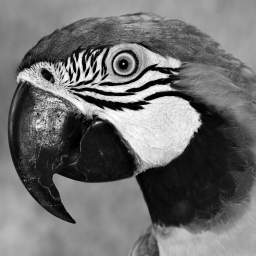}
	\endminipage
	\caption{Images used as test set for the numerical experiments.}
	\label{fig:test_set}
\end{figure}

\begin{table}[]
	\begin{tabular}{clcl}
		\multicolumn{1}{l}{}           & $\triangle t$ & M  & \multicolumn{1}{l}{ T} \\ \cline{2-4} 
		\multirow{5}{*}{$\sigma = 10$} & 0.05          & 10 & 0.5                   \\ \cline{2-4} 
		& 0.05          & 15 & 0.75                  \\ \cline{2-4} 
		& 0.1           & 10  & 1                     \\ \cline{2-4} 
		& 0.125         & 10 & 1.25                  \\ \cline{2-4} 
		& 0.125         & 12 & 1.5                   
			\end{tabular}
		\begin{tabular}{clcl}
		\multicolumn{1}{l}{}           & $\triangle t$ & M  & \multicolumn{1}{l}{ T} \\\cline{2-4} 
		\multirow{5}{*}{$\sigma = 20$} & 0.1           & 10 & 1                     \\ \cline{2-4} 
		& 0.1           & 15 & 1.5                   \\ \cline{2-4} 
		& 0.1           & 20 & 2                     \\ \cline{2-4} 
		& 0.125         & 20 & 2.5                   \\ \cline{2-4} 
		& 0.2           & 30 & 3                    
	\end{tabular}
	\caption{Time-step $\triangle t$ and number of iterations $M$ fixed for the numerical experiments. Different choices lead to  different stopping times.}
	\label{tab:parametros_fixos}
\end{table}

\begin{figure}[!htb]
\center
	\minipage{0.8\textwidth}
	\includegraphics[width=\linewidth]{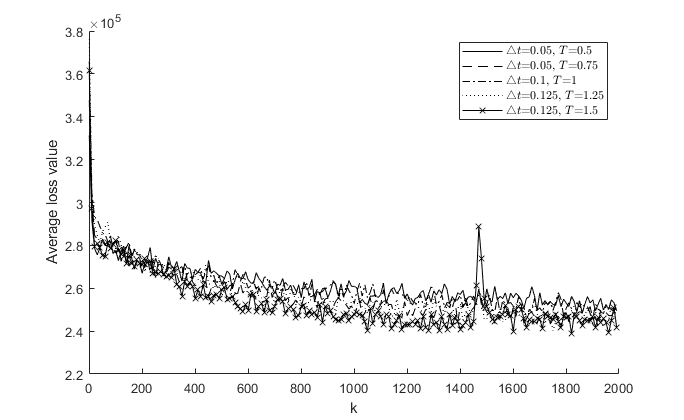}
	\endminipage
	\caption{Loss value (\ref{loss_function}) in training for gaussian denoising with $\sigma=10$, $K_{\text{max}}=2000$ and different stopping times. An initial slope gives rise to a steady decline that seems to converge to an asymptotic value shared by all stopping times.}
	\label{fig:losses_g10}
\end{figure}

\begin{figure}[!htb]
\center
	\minipage{0.45\textwidth}
	\includegraphics[width=\linewidth]{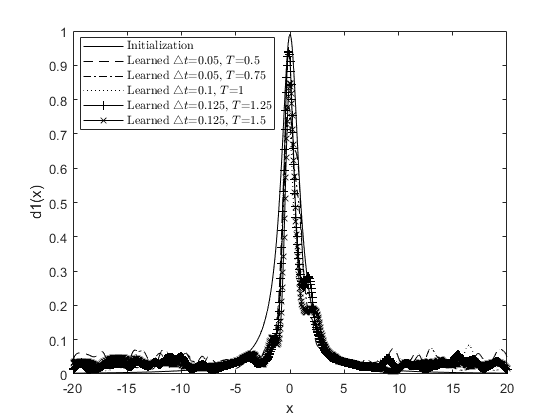}
	\endminipage
	\minipage{0.45\textwidth}
	\includegraphics[width=\linewidth]{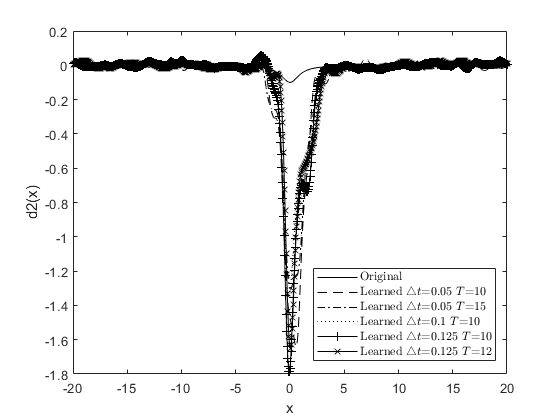}
	\endminipage \\
	
	\minipage{0.45\textwidth}
	\includegraphics[width=\linewidth]{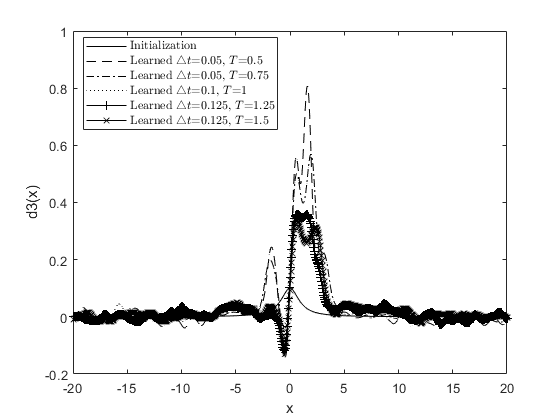}
	\endminipage
	\minipage{0.45\textwidth}
	\includegraphics[width=\linewidth]{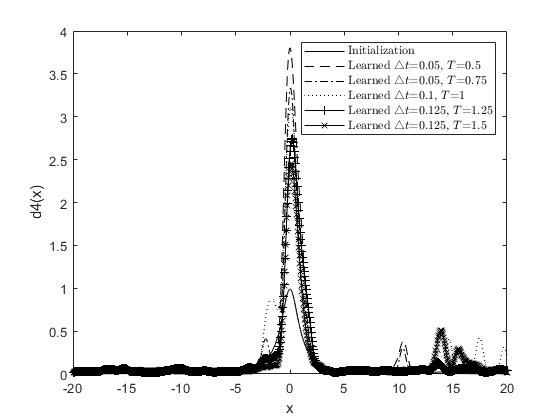}
	\endminipage
	\caption{Influence functions learning results for gaussian denoising with $\sigma=10$ and different stopping times. From left to right, top to bottom: $d_1$, $d_2$, $d_3$ and $d_4$}
	\label{fig:functions_g10}
\end{figure}

\begin{figure}[!htb]
\center
	\minipage{0.45\textwidth}
	\includegraphics[width=\linewidth]{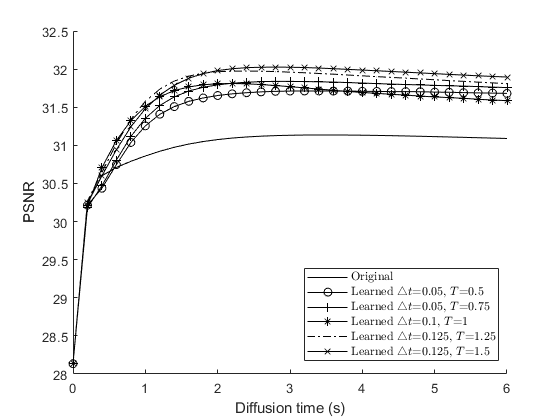}
	\endminipage
	\minipage{0.45\textwidth}
	\includegraphics[width=\linewidth]{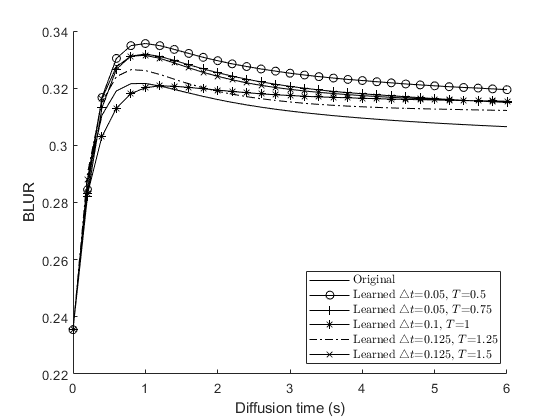}
	\endminipage
	\caption{Average of the values of PSNR (left) and Blur (right) obtained over the test set of the learning results for gaussian denoising with $\sigma=10$ and different stopping times. We obtain a considerable improvement in PSNR through training and no significant alterations in Blur.}
	\label{fig:metrics_training_set_g10}
\end{figure}

Having previously established a learning framework for the optimization of a cross-diffusion process, we now test it with some examples. We set $B=50$ and store 50 gray-scaled images (Figure \ref{fig:training_set}) that will serve as training set. The crop size for the batch extraction was fixed for all tests as $N_1=N_2=100$ and the influence functions were initialized as in the non-linear complex diffusion case \cite{AraBarSer12,AraBarSer15,GilSocZee04},
\begin{equation*}
	d_1(x)=0.99g(x), \;
	d_2(x)=-0.1g(x),  \;
	d_3(x)=0.1g(x), \;
	d_4(x)=0.99g(x), 
\end{equation*}
with 
$$g(x)=\frac{1}{1+x^2}.$$
We experimented with different values of $\triangle t$ and $T$ to see the effect of changing stopping times in this learning framework. The values chosen are summarized in (Table \ref{tab:parametros_fixos}). For the RBF interpolation (\ref{rbf}), we chose $A=[-20,20]$, $P=151$ and $\nu=0.2$. For the augmented Lagrangian parameters, we fixed for all tests $\bar{\mu}=2$, $\rho_1=6\times 10^5$, $\tau=0.5$ and $\gamma=2$. Finally, we set $K_\text{max}=2000$ for the number of  iterations of the minimization problem.

In Figure \ref{fig:losses_g10} we average the values of (\ref{loss_function}) over groups of $10$ algorithm updates for all the training procedures carried out with $\sigma=10$. The results show a similar pattern for all combination of $\triangle t$ and $M$. The algorithm reaches a significant lower loss value relatively fast (in $50$ to $100$ updates) and proceeds to improve (although at a slower pace) throughout the rest of the optimization procedure.

\begin{figure}[!htb]
\center
	\minipage{0.45\textwidth}
	\includegraphics[width=\linewidth]{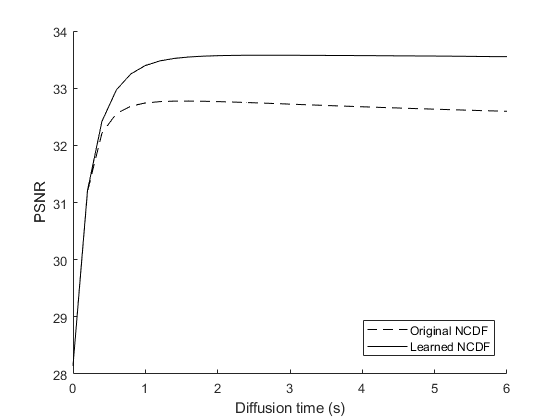}
	\endminipage
	\minipage{0.45\textwidth}
	\includegraphics[width=\linewidth]{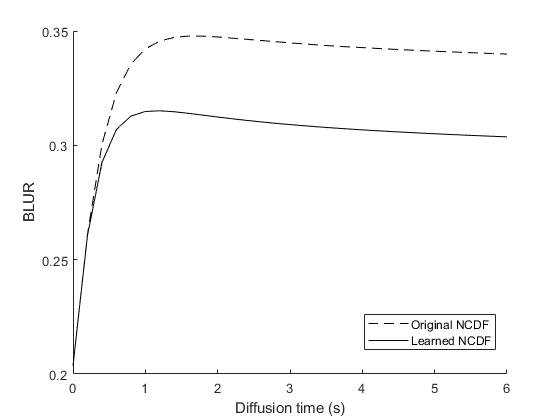}
	\endminipage \\
	
	\minipage{0.33\textwidth}
	\includegraphics[width=\linewidth]{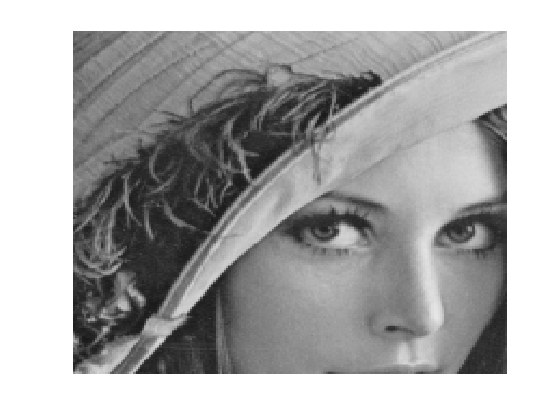}
	\vspace{-1cm}
	\caption*{\tiny Original}
	\endminipage
	\minipage{0.33\textwidth}
	\includegraphics[width=\linewidth]{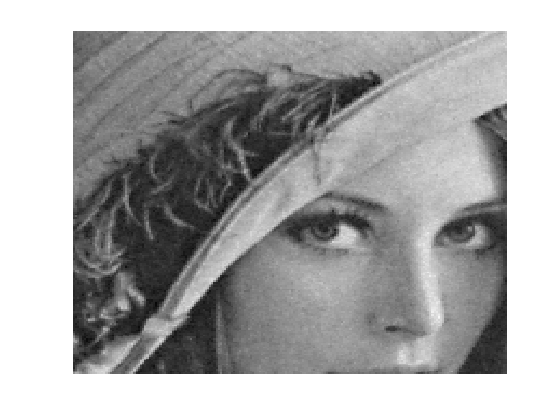}
		\vspace{-1cm}
	\caption*{\tiny NCDF}
	\endminipage
	\minipage{0.33\textwidth}
	\includegraphics[width=\linewidth]{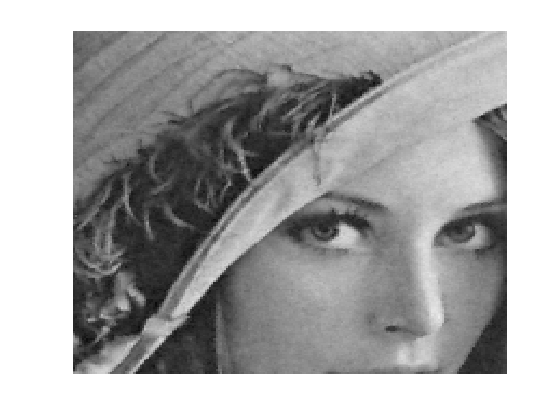}
		\vspace{-1cm}
	\caption*{\tiny Learned}
	\endminipage
	\caption{Learned NCDF results for gaussian denoising with $\sigma=10$. Top left: PSNR values of Lena denoising for original and learned NCDF. Top right: Blur values of Lena denoising for original and learned NCDF. Bottom: Lena detail of original image (left), three seconds diffusion via original NCDF (middle) and three seconds diffusion via learned NCDF (right). Considerable improvements in both metrics.}
	\label{fig:lena_g10}
\end{figure}

\begin{figure}[!htb]
\center
	\minipage{0.45\textwidth}
	\includegraphics[width=\linewidth]{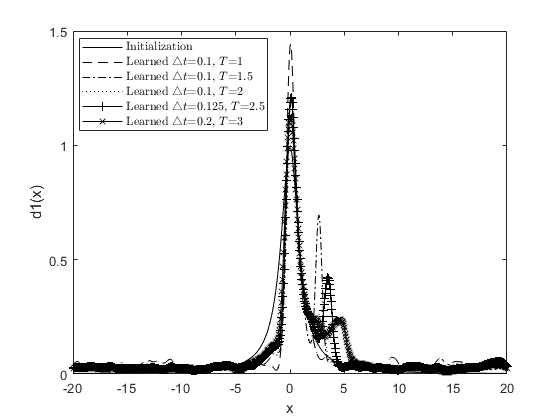}
	\endminipage
	\minipage{0.45\textwidth}
	\includegraphics[width=\linewidth]{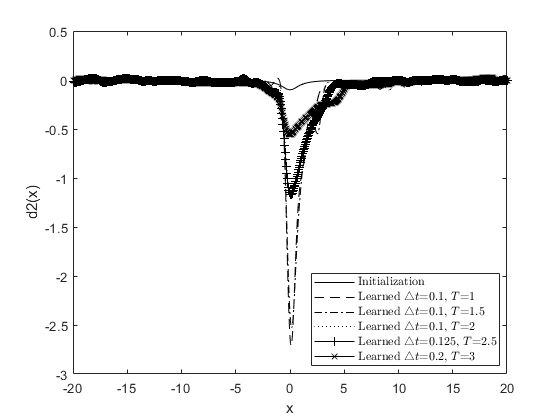}
	\endminipage \\
	
	\minipage{0.45\textwidth}
	\includegraphics[width=\linewidth]{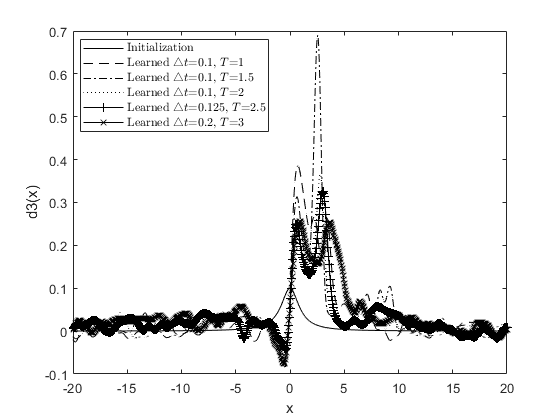}
	\endminipage
	\minipage{0.45\textwidth}
	\includegraphics[width=\linewidth]{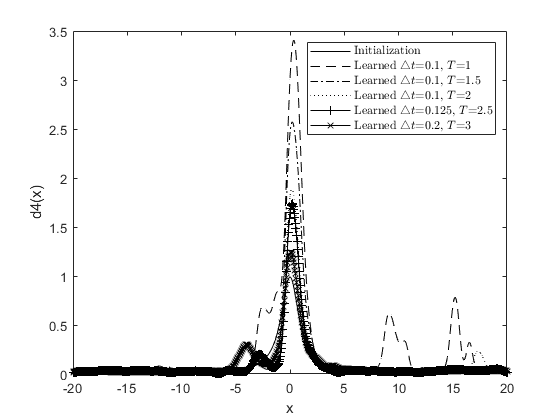}
	\endminipage
	\caption{Influence functions learning results for gaussian denoising with $\sigma=20$ and different stopping times. From left to right, top to bottom: $d_1$, $d_2$, $d_3$ and $d_4$}
	\label{fig:functions_g20}
\end{figure}

\begin{figure}[!htb]
\center
	\minipage{0.45\textwidth}
	\includegraphics[width=\linewidth]{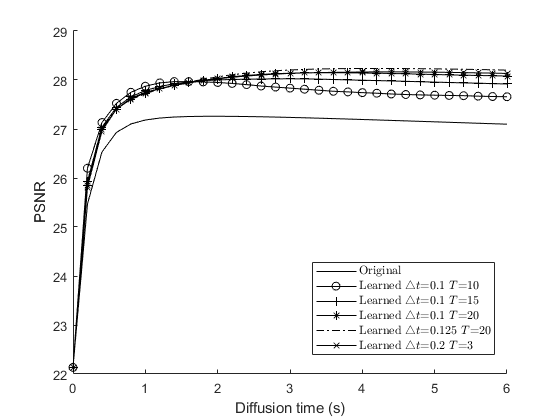}
	\endminipage
	\minipage{0.45\textwidth}
	\includegraphics[width=\linewidth]{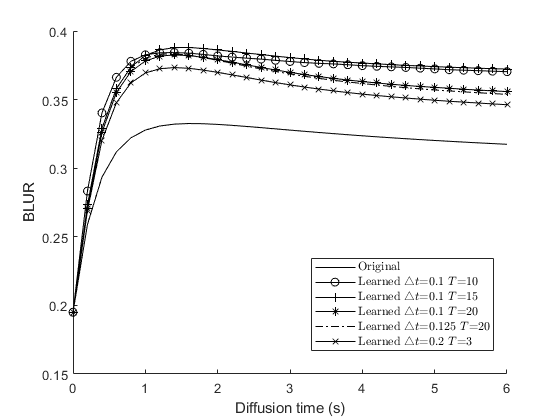}
	\endminipage
	\caption{Average of the values of PSNR (left) and Blur (right) obtained over the test set of the learning results for gaussian denoising with $\sigma=20$ and different stopping times. We obtain an improvement in PSNR through training and a slight increase in Blur.}
	\label{fig:metrics_training_set_g20}
\end{figure}

\begin{figure}[!htb]
\center
	\minipage{0.45\textwidth}
	\includegraphics[width=\linewidth]{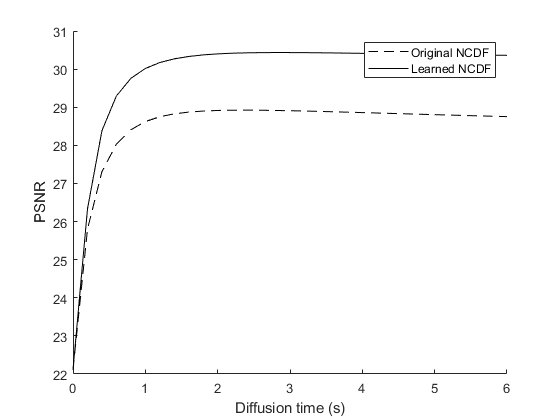}
	\endminipage
	\minipage{0.45\textwidth}
	\includegraphics[width=\linewidth]{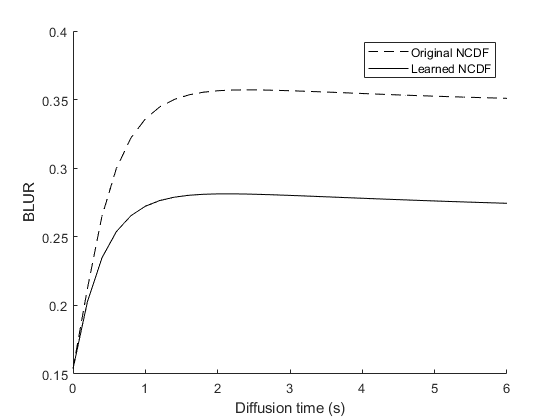}
	\endminipage \\
	
	\minipage{0.33\textwidth}
	\includegraphics[width=\linewidth]{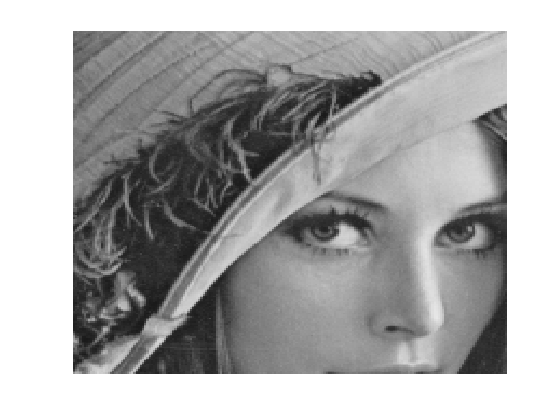}
	\vspace{-1cm}
	\caption*{\tiny Original}
	\endminipage
	\minipage{0.33\textwidth}
	\includegraphics[width=\linewidth]{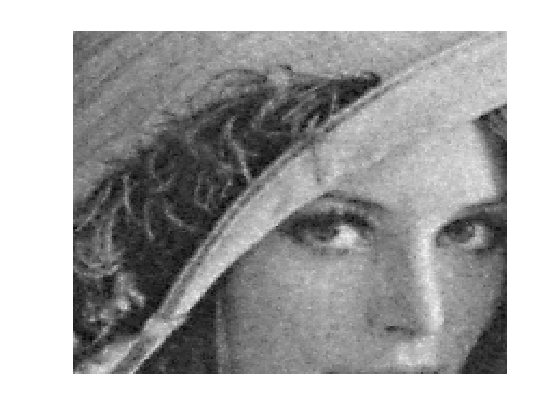}
	\vspace{-1cm}
	\caption*{\tiny NCDF}
	\endminipage
	\minipage{0.33\textwidth}
	\includegraphics[width=\linewidth]{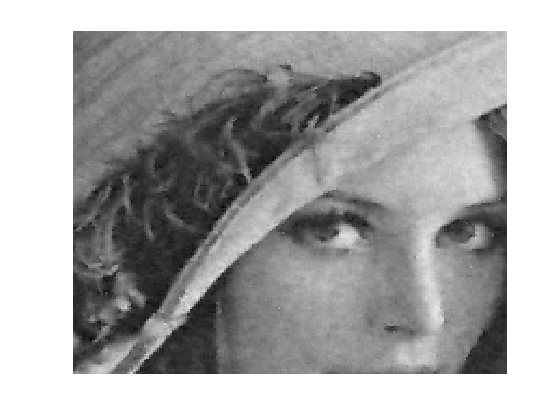}
	\vspace{-1cm}
	\caption*{\tiny Learned}
	\endminipage
	\caption{Learned NCDF results for gaussian denoising with $\sigma=20$. Top left: PSNR values of Lena denoising for original and learned NCDF. Top right: Blur values of Lena denoising for original and learned NCDF. Bottom: Lena detail of original image (left), three seconds diffusion via original NCDF (middle) and three seconds diffusion via learned NCDF (right). Considerable improvements in both metrics.}
	\label{fig:lena_g20}
\end{figure}

We describe the learned influence functions in figures \ref{fig:functions_g10} and \ref{fig:functions_g20} for $\sigma=10$ and $\sigma=20$, respectively. We notice that the learning translates into the determination of the appropriate scale of the original function in the cases of influence functions $d_1$, $d_2$ and $d_4$, and into the radical reshaping of $d_3$. The results are consistent across the different stopping times and both levels of noise.

In figures \ref{fig:metrics_training_set_g10} and \ref{fig:metrics_training_set_g20} we test the learned parameters against the non-linear complex diffusion case and with the value of $\lambda$ that provides the best PSNR in the training set. The comparisons are performed in the test set (figure $\ref{fig:test_set}$), which is made by images that were not integrated in the training set. We see in both cases a significant increase in PSNR and a slight increase in blur (\cite{blur07}). The combinations of $\triangle t$ and $M$ that achieve better performances both for $\sigma=10$ and $\sigma=20$ are those that combine into a greater value of $T$, that is, a later stopping time. The reason for this is the introduction of the reaction term, which shifts the steady-state solution into a non-trivial one (when compared to diffusion filters without such reaction), and consequently the denoising process benefits from a longer and more controlled diffusion.

Finally, we compare the action of the nonlinear complex filter against the best  learned parameters (that is, the ones trained with a larger $T$) on the widely used Lena image. The results are in figures \ref{fig:lena_g10} and \ref{fig:lena_g20} for $\sigma=10$ and $\sigma=20$, respectively. In both cases, we achieved significantly improved PSNR and blur values.

\section{Conclusions}
We have successfully  adapted a cross-diffusion model for image restoration into a learning framework, obtaining a way to automatically parametrize the cross-diffusion matrix for different images and levels of noise. 
The numerical experiments show a significant improvement of the results, measured with the PSNR and Blur metrics, obtained with the proposed learning algorithm when compared with the initialization.
By making the parallelism between neural networks and the parametrization of PDEs, we believe that this work can be transferred to a broad scope of related problems.


\section*{Acknowledgements}
This work was partially supported by the Centre for Mathematics of the University of Coimbra -- UID/MAT/00324/2019, funded by the Portuguese Government through FCT/MEC and co-funded by the European Regional Development Fund through the Partnership Agreement PT2020. The second author was supported by the FCT grant PD/BD/142956/2018.

\end{document}